\documentclass[12pt]{article}

\usepackage{amssymb}
\usepackage[leqno]{amsmath}
\usepackage{amsthm}

\usepackage{graphicx,color}

\usepackage{eso-pic}
\definecolor{lightgray}{gray}{.85}
\usepackage{subfigure}

\usepackage{xr-hyper}

\usepackage[bookmarks,pdfnewwindow]{hyperref}

\hypersetup{
colorlinks=true,
linkcolor=black,
citecolor=black,
pdfauthor={Victor Ivrii},
pdftitle={Sharp Spectral Asymptotics for 2-dimensional Schr\"odinger operator with a strong magnetic field. Note about forgotten generic case},
pdfsubject={Sharp Spectral Asymptotics},
pdfkeywords={Microlocal Analysis, Magnetic Schr\"odinger Operator},
baseurl={http://www.math.toronto.edu/ivrii/Research/preprints/},
}

\newcommand\MW{{\rm {MW}}}
\newcommand\corr{{\rm {corr}}}

\newcommand\Hess{\operatorname{Hess}}

\newcommand\bR{{\mathbb R}}

\newcommand\bZ{{\mathbb Z}}
\newcommand\cA{{\mathcal A}}

\newcommand\cE{{\mathcal E}}

\newcommand\cL{{\mathcal L}}

\newtheorem{theorem}{Theorem}

\theoremstyle{definition}

\newtheorem{remark}[theorem]{Remark}

\begin{document}

\title{%
Sharp Spectral Asymptotics for 2-dimensional Schr\"odinger operator with a strong magnetic field. Note about forgotten generic case}

\author{%
Victor Ivrii\footnote{Work was partially supported by NSERC grant OGP0138277.}
}
\date{March 27, 2005}
\maketitle

{\abstract
I consider magnetic Schr\"odinger operator in dimension $d=2$ assuming that coefficients are smooth and magnetic field is non-degenerating. Then I extend the remainder estimate $O(\mu^{-1}h^{-1}+1)$ derived in \cite{Ivr1} for the case when $V/F$ has no stationary points to the case when it has non-degenerating stationary points. If some of them are saddles and $\mu^3h\ge 2$ then asymptotics contains correction terms of magnitude 
$\mu^{-1}h^{-1}|\log \mu^3 h|$.
\endabstract}


\section{Introduction}

I consider spectral asymptotics of the magnetic Schr\"odinger operator
\begin{equation}
A= {\frac 1 2}\Bigl(\sum_{j,k}P_jg^{jk}(x)P_k -V\Bigr),\qquad P_j=D_j-\mu V_j
\label{1}
\end{equation}
where $g^{jk}$, $V_j$, $V$ are smooth real-valued functions of $x\in \bR^2$ and
$(g^{jk})$ is positive-definite matrix, $0<h\ll 1$ is a Planck parameter and 
$\mu \gg1$ is a coupling parameter. I assume that $A$ is a self-adjoint operator and all the conditions are satisfied in the ball $B(0,1)$.

In contrast to my recent papers \cite{IRO3, IRO4, IRO5} I assume that all the coefficients are very smooth; in contrast to \cite{IRO4} I consider only two-dimensional case here and in contrast to \cite{IRO6} I assume that magnetic field is non-degenerate. So I am completely in frames of section 6 \cite{Ivr1} where I just forgot to consider the case of $V/F$ having non-degenerating stationary points. My analysis will be sketchy, more details I will publish in the future. Thus this note together with Chapter 6 of \cite{Ivr1} and with \cite{IRO6} completely covers generic 2-dimensional smooth case. One can generalize these results to non-smooth case using approach of \cite{IRO3}.

Let $g=\det (g^{jk})^{-1}$, $F_{12}=\partial_{x_1}V_2-\partial_{x_2}V_1$ and
$F=|F_{12}g^{-{\frac 1 2}}|$ which is a scalar intensity of the magnetic field,
$g=\det (g^{jk})^{-{\frac 1 2}}$. I assume  that both $V$ and $F$ are disjoint from 0:
\begin{align}
&\sum_{jk} g^{jk}\xi_j\xi_k\ge \epsilon |\xi|^2\qquad\forall \xi\in \bR^2,\label{2}\\
&V\ge \epsilon_0,\label{3}\\
&F\ge \epsilon_0.\label{4}
\end{align}

In this note I am going to consider the case when $V/F$ has non-degenerate 
critical points and I will recover the same asymptotics and remainder estimate as either $\mu\le Ch^{-{\frac 1 3}}$ or $V/F$ has no saddle points in the domain in question and there will be correction terms of magnitude 
$\mu^{-1}h^{-1}|\log (\mu^3 h)|$ associated with saddle points as 
$\mu\ge 2 h^{-{\frac 1 3}}$.

I am interested in asymptotics of $\int e(x,x,0)\psi (x)\, dx$ as 
$\mu\to +\infty$, $h\to +0$ where $e(x,y,\tau)$ is the Schwartz kernel of the spectral projector of $A$ and $\psi \in C_0^\infty (B(0,{\frac 1 2})$.

\begin{theorem}\label{thm-1} Let  operator $A$ defined by $(\ref{1})$ with real-valued $g^{jk},V_j,V$ be self-adjoint in $L^2(X)$. Further  $g^{jk},V_j,V,\psi$ be smooth enough in $B(0,1)$ and conditions $(\ref{2})-(\ref{4})$ be fulfilled and there, let $B(0,1)\subset X$. Finally, let all critical points of $V/F$ in $B(0,1)$ be non-degenerate. Then

\smallskip
\noindent
(i) As $1\le \mu\le h^{-{\frac 1 3}}$ the standard asymptotics holds (i.e. $(\ref{5})-(\ref{6})$ without correction terms);

\smallskip
\noindent
(ii) As $h^{-{\frac 1 3}}\le \mu \le Ch^{-1}$
the following asymptotics holds
\begin{equation}
|\int \Bigl(e(x,x,0)-\cE^\MW (x,0)\Bigr)\psi(x)\,dx -
\sum_j \cE^\MW_\corr (x_j)\psi (x_j)|\le C\mu^{-1}h^{-1}+C
\label{5}
\end{equation}
with summation over all saddle points $x_j$ of $V/F$  where
\begin{equation}
\cE ^\MW (x,0)= {\frac 1 {2\pi}} \sum_{n\ge 0}
\theta \Bigl(\tau -V(x)-(2n+1)F\mu h \Bigr) F\mu h ^{-1}
\label{6}
\end{equation}
is    magnetic Weyl expression, and
\begin{equation}
\cE^\MW_\corr=
\varkappa   \log \Bigl(\bigl(\sigma  +\mu^{-2}\bigr)\bigl(1+\mu^{-1}h^{-1}\bigr)\Bigr)
\label{7}
\end{equation}
where 
\begin{equation}
\sigma =\sigma(x)=\min_{n\in \bZ^+} |V+(2n+1)F\mu h| 
\label{8}
\end{equation}
and $\varkappa$ is defined by $(\ref{13})$; further, as 
$C(h |\log h|)^{-1}\le \mu \le \epsilon h^{-1}$ one must include in $\cE^\MW_\corr$ 
\begin{equation}
\cE^\MW_{\corr\, 2} = 
\varkappa_2 \mu h \log \Bigl((\sigma + h^2)(1+\mu^{-1}h^{-1})\Bigr)
\label{9}
\end{equation}
again associated with saddle points.
\end{theorem}

\begin{theorem}\label{thm-2} Let  operator $A$ defined by $(\ref{1})$ with real-valued $g^{jk},V_j,V$ be self-adjoint in $L^2(X)$. Further  $g^{jk},V_j,V,\psi$ be smooth enough in $B(0,1)$ and conditions $(\ref{2}),(\ref{4})$ be fulfilled and there, let $B(0,1)\subset X$. 
Further, let $\epsilon h^{-1}\le \mu$ and $V=-(2{\bar n}+1)\mu h F+ W$ with smooth bounded $W$. Finally, let each critical point of $W/F$ in $B(0,1)$ be either non-degenerate or satisfy $|W|\ge \epsilon_0$. Then asymptotics $(\ref{5})$ holds with extra correction term $\mu h\int \varsigma \psi(x)\,dx $ 
as $\mu \le C h^{-3}|\log h|^{-1}$; for larger $\mu$ correction term contains also more complicated $O(\mu h^3|\log h|)$ terms.

\end{theorem}

\begin{remark}\label{rem-3} 
 One can drop condition (\ref{3})  by rescaling arguments after main theorem \ref{thm-1} is established.
\end{remark}

\section{Ideas of the proof: weak magnetic field case}

As $\mu \le h^{-1+\delta}$\footnote{\label{foot-1} Where here and below $\delta,\delta',\dots$ denote arbitrarily small positive  exponents.} in zone 
$\{|\nabla V|\ge \rho= C (\mu h)^{\frac 1 2}h^{-\delta}\}$ one can apply weak magnetic field approach (see section 6.3 of \cite{Ivr1}) and derive remainder estimate $O(\mu^{-1}h^{-1} + \rho^2 \mu h^{-1})$; furthermore, with logarithmic uncertainty principle replacing the standard microlocal uncertainty principle (see \cite{IRO1,IRO3}) one can derive this remainder estimate with 
$\rho = C(\mu h)^{\frac 1 2}|\log h|$. This leads to the proof of the standard asymptotics with the remainder estimate $O(\mu^{-1}h^{-1})$ as $\mu \le C(h|\log h|)^{-{\frac 1 3}}$.

Furthermore, based on the canonical form (\ref{10}) (see next section) one can prove the same asymptotics and the remainder estimate with 
$\rho= C(\mu h)^{\frac 1 2}$ and therefore achieve remainder estimate $O(\mu^{-1}h^{-1})$  as 
$\mu \le Ch^{-{\frac 1 3}}$, thus proving Theorem \ref{thm-1}(i).

\section{Ideas of the proof: intermediate and strong magnetic field cases}

To prove Theorem \ref{thm-1}(ii) and calculate correction term let me remind that according to section 6.4 of \cite{Ivr1} one can reduce microlocally operator (\ref{1}) to the canonical form
\begin{equation}
\sim \sum _{m,l,k: m+l \ge 1} a_{mnk} (x_2,\hslash D_2) 
\bigl(h^2D_1^2+\mu^2x_1^2\bigr)^m \mu^{2-2m-2l}(\mu^{-1}h)^{2k} ,\qquad \hslash=\mu^{-1}h.
\label{10}
\end{equation}
Then replacing harmonic oscillator $\bigl(h^2D_1^2+\mu^2x_1^2\bigr)$ by its eigenvalues $(2n+1)\mu h$ ($n\in\bZ^+$) one arrives to the family of  1-dimensional $\hslash$-pdos $\cA_n(x_2,\hslash D_2;\mu^{-2},\hslash)$ with symbols which modulo $O\bigl(\mu^{-2}+\mu^{-1}h\bigr)$ are 
$\Bigl(V+(2n+1)F\mu h\Bigr)\circ \Psi$ where $\Psi:\bR^2\to \bR^2$ is a map with
$|\det D\Psi|=F^{-1}$.

Since I  am interested in the energy level 0, I am most interested in the operator $\cA_n$ which is not elliptic in the point in question i.e. in operator with $n={\bar n}$ delivering minimum to $|V+(2n+1)\mu hF|$ (which I have already denoted by $\sigma$).

Furthermore, according to formula (6.6.24) of \cite{Ivr1} symbol of $\cA_n$ with $n={\bar n}$  is equal modulo $O(\mu^{-4}+h^2)$ to
\begin{align} 
&F \Bigl(-(VF^{-1}) +(2n+1)\mu h  +\mu^{-2}\omega_1\Bigr)\circ \Psi ,\label{11}\\
&\omega_1 = {\frac 1 8}\kappa V^2F^{-2} - {\frac 1 4}VF^{-1}\cL (VF^{-1})
\label{12}
\end{align}
where $\kappa$ and $\cL$ are scalar curvature and the Laplace-Beltrami operator associated with the metric $F^{-1}g^{jk}$. 

Then according to the theory of 1-dimensional operators the standard Weyl spectral asymptotics holds for each of them with the remainder estimate $O(1)$
and thus the remainder estimate for the original problem is 
$O\bigl(\mu^{-1}h^{-1}\bigr)$; however the principal part of such asymptotics 
includes the full symbol of operator, including terms of magnitude $\mu^{-2}$ and $h^{-2}$; however as $\mu \ge C h^{-{\frac 1 3}}$ one can skip terms $O(\mu^{-4})$ and $O(\mu^{-2}h^2)$ in $\cA_n$ without penalty; further, as 
$\mu \le C(h|\log h|)^{-1}$ one can skip terms $O(h^2)$ in $\cA_n$ without penalty as well.

However to preserve remainder estimate one must compensate skipping $O(\mu^{-2})$ terms in $\cA_n$  by the corresponding correction term and one can see easily that this correction term is  equal to $\kappa_0\mu^{-2}h^{-2}$ plus the correction term associated with 1-dimensional operator 
\begin{equation}
x_2 \hslash D_2 +k^{-1}(w+\mu^{-2}\omega_1)
\label{13}
\end{equation}
in zone $\{|x_2|+|\xi_2|\le \rho= C(\mu h)^{\frac 1 2}\}$ where
\begin{equation*}
k= |\det\Hess (V/F)|^{\frac 1 2},\quad
w = \bigl(-{\frac V F} + (2{\bar n}+1)\mu h\bigr),\quad\sigma=|w|
\end{equation*}
and $k, w, \sigma, \omega_1$ are calculated in the critical point in question; this latter correction term is $O\bigl(\mu^{-1}h^{-1}|\log \mu^3 h|\bigr)$ for saddle points and 
$O\bigl(\mu^{-1}h^{-1}\bigr)$ for maxima and minima and therefore only  saddle points should be considered (i.e. critical points with $\det\Hess (V/F)<0$). 

Since this asymptotics should be consistent with one obtained by weak magnetic field approach $\varkappa_0=0$ 
and the correction term in question is associated with perturbation $\mu^{-2}k^{-1}\omega_1$ in zone $|x_2|+|\xi_2|\le \rho\}$ and thus modulo $O(\mu^{-1}h^{-1})$ it is
\begin{equation*}
(2\pi) ^{-1}\mu h^{-1}F \sqrt g \times \omega_1 k^{-1} \mu^{-2}\times 
\log \Bigl({\frac \rho {|w|^{\frac 1 2}+\mu^{-1}}}\Bigr)
\end{equation*}
which can be rewritten in (\ref{7}) with 
\begin{equation}
\varkappa=
-(4\pi)^{-1}\Bigl( {\frac 1 8}\kappa V^2F^{-1} -{\frac 1 4}V\cL (VF^{-1})\Bigr) |\det\Hess (V/F) |^{-{\frac 1 2}}\sqrt g   
\label{14}
\end{equation}
calculated at this point.

Actually, this is correct only as $\mu \le C(h |\log h|)^{-1}$; 
for $C(h |\log h|)^{-1}\le \mu\le Ch^{-1}$ one should not discard an extra term 
$\omega h^2$ in $\cA_n$ but this term will contribute above $O(\mu^{-1}h^{-1})$ only as $n={\bar n}$ and it generates $\cE^\MW_{\corr\,2}$.
This leads to the proof of Theorem \ref{thm-1}(ii).

\section{Ideas of the proof:  superstrong magnetic field case}

As $\mu \ge \epsilon h^{-1}$ the same approach works but now only a single $n={\bar n}$
produces non-trivial contribution while contribution of every $n< {\bar n}$ is
$(2\pi)^{-1} \mu h^{-1}\int F\psi \,dx$ and contribution of every $n> {\bar n}$ is 0 (modulo negligible terms). So one should just repeat the same analysis where now $\rho =\epsilon $. One should not discard $\omega h^2$ in $\cA_n$ even if there are no critical points and this term produces extra correction term.
This leads to the proof of Theorem \ref{thm-2}.

\bibliographystyle{amsalpha}

\providecommand{\bysame}{\leavevmode\hbox to3em{\hrulefill}\thinspace}

\vglue .06truein

\hfill\hfill {\sl   \today \/}

\vglue .06truein

\begin{tabular}{rrl}
&{\hskip 220 pt} &Department of Mathematics,\cr
&&University of Toronto,\cr
&&100, St.George Str.,\cr
&&Toronto, Ontario M5S 3G3\cr
&&Canada\cr
&&ivrii@math.toronto.edu\cr
&&Fax: (416)978-4107\cr
\end{tabular}

\end{document}